\documentclass[amstex,12pt, amssymb]{article}

\usepackage{mathtext}
\usepackage[cp1251]{inputenc}
\usepackage[T2A]{fontenc}
\usepackage[dvips]{graphicx}
\usepackage{amsmath}
\usepackage{amssymb}
\usepackage{amsxtra}
\usepackage{latexsym}
\usepackage{ifthen}

\newcounter{lemma}[section]

\newcounter{corollary}[section]

\newcounter{remark}[section]

\newcounter{theorem}[section]

\newcounter{proposition}[section]

\newcounter{example}

\numberwithin{equation}{section}

\begin{document}

\markboth{\centerline{E.~SEVOST'YANOV, O.~DOVHOPIATYI, N.~ILKEVYCH,
V.~KALENSKA}}{\centerline{ON BOUNDARY EXTENSION OF MAPPINGS...}}

\def\cc{\setcounter{equation}{0}
\setcounter{figure}{0}\setcounter{table}{0}}

\overfullrule=0pt


\author{E.~SEVOST'YANOV, O.~DOVHOPIATYI, N.~ILKEVYCH,
V.~KALENSKA}

\title{
{\bf ON BOUNDARY EXTENSION OF MAPPINGS OF RIEMANNIAN SURFACES IN
TERMS OF PRIME ENDS}}

\date{\today}
\maketitle

\begin{abstract}
We investigate non-homeomorphic mappings of Riemannian surfaces of
Sobolev class. We have obtained some estimates of distortion of
moduli of families of curves. We have proved that, under some
conditions, these mappings have a continuous extension to a boundary
of a domain in terms of prime ends.
\end{abstract}

\bigskip
{\bf 2010 Mathematics Subject Classification: Primary 30C65;
Secondary 31A15, 31B25}

\section{Introduction}

Some important results concerning the boundary behavior of Sobolev
homeomorphisms between Riemannian surfaces were obtained
in~\cite{RV$_1$} and ~\cite{RV$_2$}. In particular, in~\cite{RV$_1$}
the authors considered the case when the domains under consideration
are locally connected at their boundary, while the
paper~\cite{RV$_2$} refers to domains of a more complex structure.
In the latter case, mappings, as a rule, do not have a pointwise
continuous boundary extension. However, the construction of prime
ends, introduced by Caratheodory, allows us to interpret this
extension in another (more successful) sense.

\medskip
In this article, we intend to abandon the condition of the
injectivity of mappings, which significantly distinguishes it
from~\cite{RV$_1$} and ~\cite{RV$_2$}. We will show that similar
classes of open-closed discrete maps also have a continuous boundary
extension. Definitions and notions used below and not mentioned in
the text, may be found in~\cite{RV$_1$}, \cite{RV$_2$}
or~\cite{Sev$_1$}.

In what follows, unless otherwise specified, the Riemannian surfaces
${\Bbb S}$ and ${\Bbb S}_*$ have hyperbolic type. In the following,
$ds_{\widetilde{h}}$ and $d\widetilde{v},$ $ds_{\widetilde{h_*}}$
and $d\widetilde{v_*}$ denote the elements of length and area on the
Riemannian surfaces ${\Bbb S}$ and ${\Bbb S}_*,$ respectively. We
also use the notation $\widetilde{h}$ for the metric on the surface
${\Bbb S},$ in particular,
$$\widetilde{B}(p_0, r):=\{p\in {\Bbb S}: \widetilde{h}(p,
p_0)<r\}, \quad\widetilde{S}(p_0, r):=\{p\in {\Bbb S}:
\widetilde{h}(p, p_0)=r\}$$
are a disk and a circle on the surface ${\Bbb S}$ centered a at a
point $p_0$ and a radius $r>0,$ respectively.

The following definitions refer to Caratheodory~\cite{Car}, see
also~\cite{KR}, \cite{GU} and earlier papers related to prime ends
(\cite{Mikl}, \cite{Su}). Recall that, a continuous mapping
$\sigma:{\Bbb I}\rightarrow {\Bbb S},$ ${\Bbb I}=(0, 1),$ is called
the {\it Jordan arc} in ${\Bbb S },$ if $\sigma(t_1)\ne\sigma(t_2)$
for $t_1\ne t_2$. Next, we will sometimes use $\sigma$ for
$\sigma({\Bbb I}),$ $\overline{\sigma}$ for $\overline{\sigma({\Bbb
I})}$ and $\partial\sigma$ for $\overline{\sigma({\Bbb
I})}\setminus\sigma({\Bbb I})$. A {\it cut} of a domain $D$ is
called either the Jordan arc $\sigma:{\Bbb I}\rightarrow D,$ ends of
which belongs to $\partial D,$ or a closed Jordan path in $D.$ A
sequence $\sigma_1,\sigma_2,\ldots,\sigma_m,\ldots$ of cuts of the
domain $D$ is called a {\it chain} if:

\medskip
(i) $\overline{\sigma_i}\cap\overline{\sigma_j}=\varnothing$ for any
$i\ne j$, $i,j= 1,2,\ldots$;

\medskip
(ii)  $\sigma_m$ {\it splits} $D$, i.e. $D\setminus \sigma_m$
consists from two components, one of which contains $\sigma_{m-1},$
and another contains $\sigma_{m+1},$

\medskip
(iii) $\widetilde{h}(\sigma_m)\rightarrow\infty$ as
$m\rightarrow\infty,$
$\widetilde{h}(\sigma_m)=\sup\limits_{p_1,p_2\in\sigma_m
}\widetilde{h}(p_1, p_2).$

\medskip
By the definition, a chain $\{\sigma_m\}$  defines the sequence of
domains $d_m\subset D$  such that $\partial\,d_m\cap
D\subset\sigma_m$ and $d_1\supset d_2\supset\ldots\supset
d_m\supset\ldots$. Two chains $\{\sigma_m\}$ and
$\{\sigma_k^{\,\prime}\}$ are called {\it equivalent}, if for any
$m=1,2,\ldots$ the domain $d_m$ contains all $d_k^{\,\prime}$
excepting a finite number and, on the other hand, for any
$k=1,2,\ldots$ the domain $d_k^{\,\prime}$ contains all $d_m$
excepting a finite number, as well. A {\it prime end} of $D$ is a
class of equivalent chains of cuts of $D$.

Let $K$ be a prime end in $D\subset {\Bbb R}^n,$ and let
$\{\sigma_m\}$ and $\{\sigma_m^{\,\prime}\}$ are chains in $K.$ In
addition,  let $d_m$ and $d_m^{\,\prime}$ are corresponding domains
with a respect to $\sigma_m$ and $\sigma_m^{\,\prime}$. Then
$$\bigcap\limits_{m=1}\limits^{\infty}\overline{d_m}\subset
\bigcap\limits_{m=1}\limits^{\infty}\overline{d_m^{\,\prime}}\subset
\bigcap\limits_{m=1}\limits^{\infty}\overline{d_m}\ ,$$ and, thus,
$$\bigcap\limits_{m=1}\limits^{\infty}\overline{d_m}=
\bigcap\limits_{m=1}\limits^{\infty}\overline{d_m^{\,\prime}}\ ,$$
i.e. the set
$$I(K)=\bigcap\limits_{m=1}\limits^{\infty}\overline{d_m}$$ depends only on
$K$ and does not depend on the chain $\{\sigma_m\}$. A set $I(K)$ is
said to be an {\it impression of a prime end} $K$. In what follows,
by $E_D$ we denote the set of all prime ends in $D,$ and
$\overline{D}_P:=D\cup E_D$ denotes the completion of $D$ by its
prime ends. Now, let us consider $\overline{D}_P$ as a topological
space in the following way. First of all, we consider that open sets
in $D$ are also open in $\overline{D}_P.$ Next, we defined a based
neighborhood of $P\subset E_D$ as a union of any domain $d,$ which
contains in some chain of $P,$ with the rest prime ends in $d.$ In
particular, in the topology mentioned above, a sequence $x_n\in D$
converges to $P\in E_D$ if and only if, for any domain $d_m,$
belonging to a sequence of domains $d_1, d_2, d_3, \ldots,$
containing in $P,$ there is $n_0=n_0(m)$ such that $x_n\in d_m$ for
$n\geqslant n_0.$

\medskip
{\it A (maximal) dilatation} of $f$ at $z$ is defined in local
coordinates by the relation
\begin{equation}\label{eq16}
K_f(z)=\frac{|f_z|+|f_{\overline{z}}|}{|f_z|-|f_{\overline{z}}|}
\end{equation}
for $J_f(z)\ne 0,$ $K_f(z)=1$ for $\Vert f^{\,\prime}(z)\Vert=0$ and
$K_f(z)=\infty$ otherwise. It is not difficult to see that, $K_f$
does not depend on local coordinates because the transition mappings
between two charts are conformal by the definition of the Riemannian
surface. The main result of the paper is the following.

\medskip
\begin{theorem}\label{th1}
{\sl\, Let $D$ and $D_*$ are domains in ${\Bbb S}$ and ${\Bbb S}_*,$
correspondingly, which have compact closures $\overline{D}\subset
{\Bbb S}$ and $\overline{D_*}\subset {\Bbb S}_*,$ while $\partial D$
and $\partial D_*$ has a finite number of components, where all
components of $\partial D_*$ are non-degenerate. Assume that,
$Q:{\Bbb S}\rightarrow (0, \infty)$ is a given function which is
measurable with a respect to the measure $\widetilde{v}$ on ${\Bbb
S},$ $Q(p)\equiv 0$ in ${\Bbb S}\setminus D.$ Let $f: D\rightarrow
D_*$ be an open, discrete and closed mapping with a finite
distortion of $D$ onto $D_*,$ such that $K_f(p)\leqslant Q(p)$ for
almost all $p\in D.$ Then $f$ has a continuous extension
$f:\overline{D}_P\rightarrow \overline{D_*}_P,$
$f(\overline{D}_P)=\overline{D_*}_P,$ if one of the following
conditions hold:

\medskip
1) the relations
\begin{equation}\label{eq45}
\int\limits_{\varepsilon}^{\varepsilon_0} \frac{dt}{\Vert
Q\Vert(t)}<\infty\,,\qquad \int\limits_0^{\varepsilon_0}
\frac{dt}{\Vert Q\Vert(t)}=\infty\,,
\end{equation}
hold at any $p_0\in \partial D,$ for some
$\varepsilon_0=\varepsilon_0(p_0)>0$ and any
$0<\varepsilon<\varepsilon_0,$ where $\Vert
Q\Vert(t):=\int\limits_{\widetilde{S}(p_0,
r)}Q(p)\,ds_{\widetilde{h}}(p)$ denotes the $L_1$-norm of the
function $Q$ over the circle $\widetilde{S}(p_0, r),$

2) the condition $Q\in FMO(\partial D)$ holds. }
\end{theorem}

\section{Preliminaries}
\setcounter{equation}{0}

In what follows, we need the following statement the proof of which
may be found in~\cite[Proposition~4.5]{Sev$_1$}, cf.
\cite[Lemma~7.4, Ch.~7]{MRSY}.

\medskip
\begin{proposition}\label{pr1A}
{\sl\, Let $p_0 \in {\Bbb S},$ let $U$ be a normal neighborhood
of$p_0,$ $0<r_1<r_2<{\rm dist}\,(p_0,
\partial U),$ let $Q(p)$ be a measurable function with a respect to
the measure $\widetilde{v},$ $Q:{\Bbb S}\rightarrow [0, \infty],$
$Q\in L^1(U).$ Set $\widetilde{A}=\widetilde{A}(p_0, r_1,
r_2)=\{p\in {\Bbb S}: r_1<\widetilde{h}(p, p_0)<r_2\},$ $\Vert
Q\Vert(r)=\int\limits_{\widetilde{S}(p_0,
r)}Q(p)\,\,ds_{\widetilde{h}}(p),$
\begin{equation}\label{eq8B}
 \eta_0(r):=\frac{1}{I\cdot\Vert Q\Vert(r)}\,,
\end{equation}
where
$$I=I(p_0,r_1,r_2):=\int\limits_{r_1}^{r_2}\
\frac{dr}{\Vert Q\Vert(r)}\,.$$
Then
\begin{equation}\label{eq10AA}
\frac{1}{I}=\int\limits_{\widetilde{A}(p_0, r_1, r_2)} Q(p)\cdot
\eta_0^2(\widetilde{h}(p, p_0))\
d\widetilde{h}(p)\leqslant\int\limits_{\widetilde{A}(p_0, r_1, r_2)}
Q(p)\cdot \eta^2(\widetilde{h}(p, p_0))\ d\widetilde{v}(p)
\end{equation}
for any Lebesgue measurable function $\eta:(r_1,r_2)\rightarrow
[0,\infty]$ such that
\begin{equation}\label{eq8AA}\int\limits_{r_1}^{r_2}\eta(r)dr=1\,.
\end{equation}
}
\end{proposition}

\medskip
Given a mapping $f:D\rightarrow {\Bbb S}_*$ and a set $E\subset
\overline{D}\subset{\Bbb S},$ we put
$$C(f, E)=\{y\in {\Bbb S}_*:\exists\, x\in E, x_k\in D: x_k\rightarrow x, f(x_k)
\rightarrow y, k\rightarrow\infty\}\,.$$
The following statement holds.

\medskip
\begin{proposition}\label{pr1}
{\sl\, Assume that, a domain $D\subset {\Bbb S}$ has a finite number
of boundary components $\Gamma_1, \Gamma_2,\ldots, \Gamma_n\subset
\partial D.$ Then:

\medskip
1) for any $\Gamma_i,$ $i=1,2,\ldots, n$ there is a neighborhood
$U_i\subset {\Bbb S}$ and a conformal mapping $H$ of $U^*_i:=U_i\cap
D$ onto $R=\{z\in {\Bbb C}: 0\leqslant r_i<|z|<1\}$ such that
$\gamma_i:=\partial U_i^{\,*}\cap D$ is a closed Jordan path
$$C(H, \gamma_i) = \{z\in {\Bbb C}: |z| = 1\};\quad  C(H, \Gamma_i) =
\{z\in {\Bbb C}: |z| = r_i\}\,,$$
while $r_i=0$ if and only if $\Gamma$ degenerates into a point.
Moreover, $H$ extends to a homeomorphism of $\overline{U_i^{\,*}}_P$
onto $\overline{R},$ see \cite[Lemma~2]{RV$_2$};

\medskip
2) a space $\overline{D}_P$ is metrizable with some metric
$\rho:\overline{D}_P\times\overline{D}_P\rightarrow {\Bbb R}$ such
that, the convergence of any sequence $x_n\in D,$ $n=1,2,\ldots ,$
to some prime end $P\in E_D$ is equivalent to the convergence $x_n$
in one of spaces $\overline{U_i^{\,*}}_P,$ see
\cite[Remark~2]{RV$_2$};

\medskip
3) any prime end $P\in E_D$ contains a chain of cuts $\sigma_m,$
$m=1,2,\ldots,$ which belong to spheres $\widetilde{S}(z_0, r_m),$
$r_m\rightarrow\infty$ as $m\rightarrow\infty,$
see~\cite[Remark~1]{RV$_2$};

\medskip
4) for any $P\subset E_D$ its impression $I(P)$ is a continuum in
$\partial D,$  while there is some unique $1\leqslant i\leqslant n$
such that $I(P)\subset \Gamma_i,$ see \cite[Proposition~1,
Remark~1]{RV$_2$}.}

\end{proposition}

\medskip
The technique for proving the main result is based on using modulus
of families of paths, the distortion of which under mappings of the
Sobolev class is known. Proceeding from this, we consider some
(wider) class of mappings for which the required distortion of
modulus is satisfied. Everywhere below, $M(\cdot)$ is the modulus of
families of paths on ${\Bbb S}$ (see, for example,
\cite{RV$_1$}--\cite{Sev$_1$}). Let $\rho\colon{\Bbb
S}\rightarrow[0, \infty]$ is a function which is measurable with
respect to the area $\widetilde{v}.$ We say that, $\rho$ is {\it
extensively admissible} for $\Gamma,$ abbr. $\rho\in{\rm ext}\,{\rm
adm}\,\Gamma,$ if the ratio
$$\int\limits_{\gamma}\rho\,ds_{\widetilde{h}}(p) \geqslant 1$$
holds for all locally rectifiable paths $\gamma\in
\Gamma\setminus\Gamma_0,$ while $M(\Gamma_0)=0.$ The following class
is a generalization of quasiconformal mappings in Gehring sense
(see, e.g., \cite[Chapter~9]{MRSY}). Let $D$ and $D_{\,*}$ be
domains in ${\Bbb S}$ and ${\Bbb S}_*,$ respectively, and let
$Q\colon D\rightarrow(0,\infty)$ be a measurable function with a
respect to $\widetilde{v}$ on ${\Bbb S}.$ We say that, $f\colon
D\rightarrow D_{\,*}$ is a {\it lower $Q$-mapping at a point}
$p_0\in\overline{D},$ if there is
$\varepsilon_0=\varepsilon_0(p_0)>0,$
$\varepsilon_0<d_0=\sup\limits_{p\in D}\widetilde{h}(p, p_0),$ such
that
\begin{equation}\label{eq1A}
M(f(\Sigma_{\varepsilon}))\geqslant \inf\limits_{\rho\in{\rm
ext\,adm}\,\Sigma_{\varepsilon}}\int\limits_{D\cap
\widetilde{A}(p_0, \varepsilon,
\varepsilon_0)}\frac{\rho^2(p)}{Q(p)}\,d\widetilde{v}(p)
\end{equation}
for any ring $\widetilde{A}(p_0, \varepsilon, \varepsilon_0)=\{p\in
{\Bbb S}: \varepsilon<\widetilde{h}(p, p_0)<\varepsilon_0\},$ where
$\Sigma_{\varepsilon}$ denotes the family of all intersections of
circles $\widetilde{S}(p_0, r)=\{p\in {\Bbb S}: \widetilde{h}(p,
p_0)=r\}$ with $D,$ $r\in (\varepsilon, \varepsilon_0).$

\medskip
In many cases, we need to verify the property~(\ref{eq1A}) without
of a verification of infinitely many inequalities. Such a
possibility follows by the following statement (cf.
\cite[Theorem~9.2]{MRSY} and \cite[Lemma~4.2]{IS}), the proof of
which may be found in~\cite[Lemma~2.3]{Sev$_1$}.

\medskip
 \begin{lemma}\label{lem4A}
{\sl\, Let $D$ and $D_{\,*}$ be domains in ${\Bbb S}$ and ${\Bbb
S}_*,$ respectively, $p_0\in\overline{D}$ and $Q\colon
D\rightarrow(0,\infty)$ be a given function. Then $f\colon
D\rightarrow D_{\,*}$ is a lower $Q$-mapping at a point $p_0$ if and
only if there is $0<d_0<\sup\limits_{p\in D}\widetilde{h}(p, p_0)$
such that
\begin{equation}\label{eq15}
M(f(\Sigma_{\varepsilon}))\geqslant
\int\limits_{\varepsilon}^{\varepsilon_0}
\frac{dr}{\Vert\,Q\Vert(r)}\quad\forall\
\varepsilon\in(0,\varepsilon_0)\,,\ \varepsilon_0\in(0,d_0)\,,
\end{equation}
where $\Sigma_{\varepsilon}$ denotes the family of all intersections
$\widetilde{S}(p_0, r)$ with a domain $D,$ $r\in (\varepsilon,
\varepsilon_0),$ in addition,
$$
\Vert Q\Vert(r)=\int\limits_{D(p_0,r)}Q(p)\,ds_{\widetilde{h}}(p)$$
denotes $L_{1}$-norm of the function $Q$ over $D\cap
\widetilde{S}(p_0,r)=D(p_0,r)=\{p\in D\,:\, \widetilde{h}(p,
p_0)=r\}$.}
\end{lemma}

\medskip
The following lemma holds.

\medskip
\begin{lemma}\label{lem1}
{\sl\, Let $D$ and $D_*$ be domains in ${\Bbb S}$ and ${\Bbb S}_*,$
correspondingly, which have compact closures $\overline{D}\subset
{\Bbb S}$ and $\overline{D_*}\subset {\Bbb S}_*,$ while $\partial D$
and $\partial D_*$ consist of finite boundary components, and all
components of $\partial D_*$ are non-degenerate. Assume that,
$f:D\rightarrow D_*,$ $f(D)=D_*,$ be a closed open discrete mapping.
Then:

\medskip
1) $C(f, P)$ is a continuum in $\partial D_*,$ where
$$C(f, P)=\{y\in {\Bbb S}_*: \exists\,x_k\in D: x_k\rightarrow P,
f(x_k)\rightarrow y, k\rightarrow\infty\}\,.$$
In particular, there is a unique component $\Gamma\subset
\partial D_*$ such that $C(f, P)\subset\Gamma;$

\medskip
2) if $P\subset E_D$ and $d_k,$ $k=1,2,\ldots,$ be a sequence of
domains corresponding to $P$ and $U\subset {\Bbb S}_*$ be a
neighborhood of $\Gamma$ from item~1) of Proposition~\ref{pr1}, then
there is $s_0\in {\Bbb N}$ such that
\begin{equation}\label{eq5B}
f(d_k)\subset U^*\quad\forall\,\, k\geqslant s_0\,,
\end{equation}
where $U^*:=U\cap D_*.$
 }
\end{lemma}

\medskip
\begin{proof}
Let us to prove that $C(f, P)$ is a continuum in $\partial D_*,$
where
$$C(f, P)=\{y\in {\Bbb S}_*: \exists\,x_k\in D: x_k\rightarrow P,
f(x_k)\rightarrow y, k\rightarrow\infty\}\,.$$ For this goal, let us
to show that
\begin{equation}\label{eq3}
C(f, P)=\bigcap\limits_{k=1}^{\infty}\overline{f(d_k)}\,,
\end{equation}
where $d_k,$ $k=1,2,\ldots ,$ is a sequence of domains of cuts
corresponding to a prime end $P.$ Indeed, let $y\in C(f, P),$ then
$y=\lim\limits_{k\rightarrow\infty}y_k,$ $y_k\rightarrow P$ as
$k\rightarrow\infty.$ We may consider that $y_k=f(x_k),$ $x_k\in
d_k.$ Then, for any $m\in {\Bbb N}$ there is $k_0=k_0(m)$ such that
$x_k\in d_m$ for $k\geqslant k_0,$ because the sequence $d_m$ is
decreasing. It follows from this that, $y\in \overline{f(d_k)}$ for
any $k=1,2,\ldots.$ Thus $C(f,
P)\subset\bigcap\limits_{k=1}^{\infty}\overline{f(d_k)}.$ On the
other hand, let $y\in
\bigcap\limits_{k=1}^{\infty}\overline{f(d_k)}.$ Then, for a given
$k\in {\Bbb N},$ we obtain that
$y=\lim\limits_{m\rightarrow\infty}y^{(k)}_m,$ where $y^{(k)}_m\in
f(d_k),$ $m=1,2,\ldots .$ Then there are $x^{(k)}_m\in d_k,$
$m=1,2,\ldots ,$ such that $f(x^{(k)}_m)\rightarrow y$ as
$m\rightarrow\infty.$ Then, for a number $1/2^{\,k},$ there is
$m=m_k\in {\Bbb N}$ such that $\widetilde{h}_*(f(x^{(k)}_{m_k}),
y)<1/2^{\,k}.$ By the definition, a sequence $x^{(k)}_{m_k}$
converges to $P$ as $k\rightarrow\infty$ and
$f(x^{(k)}_{m_k})\rightarrow y$ as $k\rightarrow\infty,$ i.e., $y\in
C(f, P).$ Thus, $C(f,
P)\subset\bigcap\limits_{k=1}^{\infty}\overline{f(d_k)},$
$\bigcap\limits_{k=1}^{\infty}\overline{f(d_k)}\subset C(f, P)$ and,
consequently, the relation~(\ref{eq3}) is established. Then,
by~\cite[Theorem~5.II.5]{Ku} $C(f, P)$ is a continuum.

\medskip
It remains to show that $C(f, P)\subset \partial D_*.$ Observe that,
$C(f, P)\ne \varnothing$ because $\overline{D_*}$ is a compactum by
the assumption.  Let $y\in C(f, P),$ then
$y=\lim\limits_{k\rightarrow\infty}y_k,$ $y_k\rightarrow P$ as
$k\rightarrow\infty$ and $y_k=f(x_k),$ $x_k\in d_k.$ Without loss of
generality, by the compactness of $\overline{D},$ we may consider
that $x_k$ converges to $x_0$ as $k\rightarrow\infty.$ Then, by
Proposition~\ref{pr1}, since $x_0\in
\bigcap\limits_{k=1}^{\infty}\overline{d_k},$ we have that $x_0\in
I(P)\subset
\partial D.$ Since $f$ is an open, discrete and closed mapping, it
is boundary preserving. Now, the sequence $f(x_k)=y_k$ may converges
only to a boundary point as $k\rightarrow\infty,$ i.e., $y\in
\partial D_*.$ The item~1) of Lemma~\ref{lem1} is established.

\medskip
Let us to prove item~2). Let $U\subset {\Bbb S}_*$ be an open of
$\Gamma$ which corresponds to Proposition~\ref{pr1}. In other words,
there is a conformal mapping $H$ of $U^*:=U\cap D_*$ onto the ring
$R=\{z\in {\Bbb C}: 0<r<|z|<1\}$ such that $\gamma:=\partial
U^{\,*}\cap D$ is a closed Jordan path,
$$C(H, \gamma) = \{z\in {\Bbb C}: |z| = 1\};\quad  C(H, \Gamma) =
\{z\in {\Bbb C}: |z| = r\}\,.$$

Let us to prove~(\ref{eq5B}). Assume the contrary. Then there is an
increasing sequence of numbers $k_l,$ $l=1,2, \ldots,$ and a
sequence $y_{k_l}\in f(d_{k_l})$ such that $y_{k_l}\in D_*\setminus
U^*$ for any $l\in {\Bbb N}.$ By the compactness of $\overline{D_*}$
we may assume that $y_{k_l}$ converges to some point $y_0$ as
$l\rightarrow\infty.$ Then $y_0\in \Gamma$ by the inclusion~$C(f,
P)\subset \Gamma.$ Let $\varepsilon_1>0$ be such that $B(y_0,
\varepsilon_1)\subset U;$ this $\varepsilon_1$ exists because $U$ is
a neighborhood of $\Gamma.$ Then $y_{k_l}\in B(y_0,
\varepsilon_1)\cap D_*\subset U^{\,*}$ for large $l\in{\Bbb N},$
that contradicts with $y_{k_l}\in D_*\setminus U^*$ for $l\in {\Bbb
N}.$ The contradiction obtained above proves the
relation~(\ref{eq5B}).~$\Box$
\end{proof}

\medskip
An analog of the following statement is proved for homeomorphisms
in~\cite[Lemma~4]{RV$_2$}, cf.~\cite[Lemma~3]{KR} and
\cite[Lemma~5.1]{GRY}.

\medskip
\begin{theorem}\label{th2} {\,\sl Let
$D$ and $D_*$ domains in ${\Bbb S}$ and ${\Bbb S}_*,$
correspondingly, which have compact closures $\overline{D}\subset
{\Bbb S}$ and $\overline{D_*}\subset {\Bbb S}_*,$ while $\partial D$
and $\partial D_*$ consist of a finite number of components, and all
components of $\partial D_*$ are non-degenerate. Assume that,
$Q:{\Bbb S}\rightarrow (0, \infty)$ is a given function which is
measurable with a respect to $\widetilde{v}$ on ${\Bbb S},$
$Q(p)\equiv 0$ in ${\Bbb S}\setminus D.$ Let $f:D\rightarrow D_*,$
$D_*=f(D),$  be a lower $Q$-mapping at any point $p_0\in
\partial D,$ and let $f$ be an open, discrete and closed. Then $f$
has a continuous extension $f:\overline{D}_P\rightarrow
\overline{D_*}_P,$ $f(\overline{D}_P)=\overline{D_*}_P,$ whenever
one of the following conditions hold:

\medskip
1) either the relation
\begin{equation}\label{eq10}
\int\limits_{\varepsilon}^{\varepsilon_0} \frac{dt}{\Vert
Q\Vert(t)}<\infty\,,\qquad \int\limits_{0}^{\varepsilon_0}
\frac{dt}{\Vert Q\Vert(t)}=\infty
\end{equation}
holds for any $p_0\in \partial D,$ for some
$\varepsilon_0=\varepsilon_0(p_0)>0$ and all
$0<\varepsilon<\varepsilon_0,$ where $\Vert
Q\Vert(t):=\int\limits_{\widetilde{S}(p_0,
r)}Q(p)\,ds_{\widetilde{h}}(p)$ denotes the $L_1$-norm of the
function $Q$ over the circle $\widetilde{S}(p_0, r),$

2) or $Q\in FMO(\partial D).$ }
\end{theorem}

\medskip
\begin{proof}
Let us firstly prove that $f$ has a continuous extension
$f:\overline{D}_P\rightarrow \overline{D_*}_P.$ Let us consider the
case~1), i.e., when the relations~(\ref{eq10}) hold. Put $P\in E_D.$

\medskip
\textbf{1)} By the item~1) of Lemma~\ref{lem1}, the set $C(f, P)$ is
a continuum in $\partial D_*.$ Then there is a component
$\Gamma\subset
\partial D_*$ which contains $C(f, P).$ Let $U\subset {\Bbb S}_*$
be a neighborhood $\Gamma$ which corresponds to
Proposition~\ref{pr1}, and let $H$ be a corresponding conformal
mapping of a domain $U^*:=U\cap D_*$ onto the ring $R=\{z\in {\Bbb
C}: 0<r<|z|<1\}$ such that $\gamma:=\partial U^{\,*}\cap D$ is a
closed Jordan path,
$$C(H, \gamma) = \{z\in {\Bbb C}: |z| = 1\};\quad  C(H, \Gamma) =
\{z\in {\Bbb C}: |z| = r\}\,.$$
By Proposition~\ref{pr1} there is a chain of cuts $\sigma_n,$
corresponding to a prime end $P,$ which belongs to spheres
$\widetilde{S}(p_0, r_n),$ $p_0\in
\partial D,$ $r_n\rightarrow 0$ as $n\rightarrow\infty.$
Let $d_n,$ $n=1,2,\ldots ,$ be a sequence of domains corresponding
to cuts $\sigma_n.$ By the inclusion~(\ref{eq5B}) we may consider
that $f(d_1)\subset U^*.$ Now, we set $\widetilde{f}:=f|_{d_1},$
$g:=H\circ \widetilde{f},$ $g:d_1\rightarrow R,$ $g(d_1)\subset R.$
Observe that, $R$ is a domain, any point of which has a sufficiently
small neighborhood the intersection of which with $R$ is
quasiconformally equivalent to the unit disk (besides the direct
arguing, this statement may be obtained by the corresponding
V\"{a}is\"{a}l\"{a}'s result~\cite[Theorem~17.12]{Va$_2$}, since $R$
is a union of two circles which are $C^1$-manifolds. In this
context, we also mention~\cite[sect.~2.2]{Na$_3$}
and~\cite[Remark~1.5]{Na$_1$}). Then, by~\cite[Theorem~4.1]{Na$_3$}
and due to~\cite[Remark~2]{RV$_2$}, we may consider that
$\overline{R}_P=\overline{R}.$ In this case, for the proof of
Theorem, it is sufficient to establish the continuous extension
$\overline{g}:d_1\cup\{P\}\rightarrow \overline{R}.$

\medskip
\textbf{2)} Moreover, by the compactness of $\overline{R},$ it is
sufficiently to prove that the set
$$L=C(g, P):=\left\{y\in \partial R: y=\lim\limits_{m\rightarrow\infty}g(p_m),
p_m\rightarrow P, p_m\in d_1\right\}$$ consists from a unique point
$y_0\in\partial R.$ The mapping $g,$ as usual, is open and discrete
in $d_1,$ but is not necessary closed. Let us to show that, $g$
satisfies the relation
\begin{equation}\label{eq9A} M(g(\Sigma^1_{\varepsilon}))\geqslant
\int\limits_{\varepsilon}^{\varepsilon_1}
\frac{dr}{\Vert\,Q\Vert(r)}\quad\forall\
\varepsilon\in(0,\varepsilon_1)\,,\ \varepsilon_1\in(0, r_1)\,,
\end{equation}
where $\Sigma^1_{\varepsilon}$ denotes the family of all
intersections of circles $\widetilde{S}(p_0, r)=\{p\in {\Bbb S}:
\widetilde{h}(p, p_0)=r\}$ with $d_1,$ $r\in (\varepsilon,
\varepsilon_1).$
To proof this fact, let us to show that
\begin{equation}\label{eq6}
\widetilde{S}(p_0, r)\cap D> \widetilde{S}(p_0, r)\cap d_1, \quad
r<r_1\,.
\end{equation}
(Here and below the notation $\Gamma_1>\Gamma_2$ denotes that, for
any dished line $\alpha\in \Gamma_1,$
$\alpha:\bigcup\limits_{i=1}^{\infty}(a_i, b_i)\rightarrow {\Bbb
S},$ there is a dashed line $\beta\in \Gamma_2,$ where
$\beta:\bigcup\limits_{i,k=1}^{\infty}(a_{ik}, b_{ik})\rightarrow
{\Bbb S},$ $\bigcup\limits_{k=1}^{\infty}(a_{ik}, b_{ik})\subset
(a_i, b_i),$ $\alpha|_{(a_{ik}, b_{ik})}=\beta$ for any
$i=1,2,\ldots,$ $k=1,2,\ldots$ and, besides that, al least one
interval $(a_{ik}, b_{ik})$ is not empty).

\medskip
Indeed, we put $0<r<r_1.$ Then there is $i\in {\Bbb N}$ such that
$r_i<r.$ Let $\sigma_i\subset \widetilde{S}(p_0, r_i)\cap d_1$ be a
cut corresponding to a domain $d_i.$ Join any point
$\omega\in\sigma_i$ with a point $w\in \sigma_1$ in $D$ with a path
$\alpha_i,$ $|\alpha_i|\subset D.$ Without loss of generality, we
may consider that $\alpha_i$ belongs to $d_1$ instead of its
endpoint, because $\partial d_1\cap D\subset \sigma_1.$ Observe that
$|\alpha_i|\cap \widetilde{B}(p_0, r)\ne\varnothing\ne
|\alpha_i|\cap ({\Bbb S}\setminus \widetilde{B}(p_0, r)),$
therefore, by~\cite[Theorem~1.I.5, $\S\,46$]{Ku} $|\alpha_i|\cap
\widetilde{S}(p_0, r)\ne\varnothing.$ It follows from this that,
$\widetilde{S}(p_0, r)\cap d_1\ne\varnothing.$

\medskip
Let now $\alpha:=\widetilde{S}(p_0, r)\cap D$ be a dished line
$\alpha:\bigcup\limits_{i=1}^{\infty}(a_i, b_i)\rightarrow D,$ where
there at least one non-empty interval in its system $(a_i, b_i).$ By
the proved above, $\beta:=\widetilde{S}(p_0, r)\cap
d_1\ne\varnothing,$ therefore there is at most countable of
intervals $(c_k, d_k),$ $k=1,2,\ldots ,$ such that
$\beta:\bigcup\limits_{k=1}^{\infty}(c_k, d_k)\rightarrow d_1$ and
the interval $(c_k, d_k)$ is not empty at least for some $k\in{\Bbb
N}.$ By the definition, for any $k\in {\Bbb N}$ there is $i\in {\Bbb
N}$ such that $(c_k, d_k)\subset (a_i, b_i).$  Denote $(a_{ik},
b_{ik}):=(a_i, b_i)\cap (c_k, d_k),$ and observe that the interval
$(a_{ik}, b_{ik})$ is not empty  at least for some $i\in{\Bbb N}$
and $k\in {\Bbb N}.$ Observe also that,
$\bigcup\limits_{k=1}^{\infty}(a_{ik}, b_{ik})\subset (a_i, b_i)$
and $\alpha|_{(a_{ik}, b_{ik})}=\beta,$ that proves~(\ref{eq6}).

\medskip
It follows from~(\ref{eq6}) that $f(\widetilde{S}(p_0, r)\cap D)>
f(\widetilde{S}(p_0, r)\cap d_1)=\widetilde{f}(\widetilde{S}(p_0,
r)\cap d_1),$ where $\widetilde{f}:=f|_{d_1}.$ Let
$\Sigma_{\varepsilon}$ be a family of all intersections of circles
$\widetilde{S}(p_0, r)=\{p\in {\Bbb S}: \widetilde{h}(p, p_0)=r\}$
with $D.$ Then by~\cite[теорема~1(c)]{Fu}
$M(\widetilde{f}(\Sigma^1_{\varepsilon}))\geqslant
M(f(\Sigma_{\varepsilon}))$ and, consequently, by
Lemma~(\ref{lem4A})
\begin{equation}\label{eq10B}
 M(\widetilde{f}(\Sigma^1_{\varepsilon}))\geqslant
\int\limits_{\varepsilon}^{\varepsilon_1}
\frac{dr}{\Vert\,Q\Vert(r)}\quad\forall\
\varepsilon\in(0,\varepsilon_1)\,,\ \varepsilon_1\in(0, r_1)\,.
\end{equation}
In this case, (\ref{eq9A}) follows by~(\ref{eq10B}), because
$g=H\circ \widetilde{f}$ and $H$ is a conformal mapping preserving
the family of paths with a respect to Lebesgue measure on the plane
(see, e.g., \cite[Theorem~8.1]{Va$_2$}, see also the corresponding
result about equality of the moduli of families of paths in the
hyperbolic and Euclidean metrics and
measures~\cite[Remark~5.2]{Sev$_2$}. On this occasion we also
mention on~\cite[Remark~1]{VR}, where the notion of the modulus of
families of paths is given in some another (equivalent) way.

\medskip Put $\delta\in (0, r_1)$ and set $\Gamma^{\,\delta}_n:=\bigcup\limits_{r\in (r_n, \delta)}
g(\widetilde{S}(p_0, r)\cap d_1),$ where the union must be
understood not in the theoretical-set sense, but namely as a family
of paths ''from $r_n$ to $\delta$''. By~(\ref{eq10B}) and due
to~(\ref{eq10}) it follows that
\begin{equation}\label{eq40A}
M(\Gamma^{\,\delta}_n)\rightarrow\infty\,,\quad
n\rightarrow\infty\,.
\end{equation}
\textbf{3)} Let us prove by contradiction, i.e., assume that $g$ has
no a limit as $p\rightarrow P.$ Then we may find at least two
sequences $p_n, p_n^{\,\prime}\in d_n,$ $n=1,2,\ldots ,$ and two
points $y\ne y_*,$ $y, y_*\in R$ such that $g(p_n)\rightarrow y$ and
$g(p^{\,\prime}_n)\rightarrow y_*$ as $n\rightarrow\infty.$ Join the
points $p_n$ and $p^{\,\prime}_n$ by a path $\gamma_n$ in a domain
$d_n.$ Let $r_0:=|y-y_*|$ and $U_0:=B(y, r_0/2).$ Observe that, a
boundary of $R$ is strongly accessible since $R$ has a finite number
of components and is a finitely connected on the boundary (see,
e.g., \cite[Theorem~6.2 and Corollary~6.8]{Na$_2$}). Therefore, for
a neighborhood $U_0$ of $y$ there is an another neighborhood
$V\subset U_0$ of this point, a compactum $K\subset R$ and a number
$\delta>0$ such that the relation
\begin{equation}\label{eq4}
M(\Gamma(E, K, R))\geqslant\delta
\end{equation}
holds for any continuum $E\subset R,$ $E\cap\partial
U\ne\varnothing\ne E\cap\partial V.$

Let $C_n:=g(|\gamma_n|).$ Observe that, $\partial U\cap
|C_n|\ne\varnothing$ for sufficiently large $n\in {\Bbb N},$
see~\cite[Theorem~1.I.5, $\S\,46$]{Ku} (as usually, $|C_n|$ denotes
the locus of a path $C_n$). Then, by the condition~(\ref{eq4}) we
obtain that
\begin{equation}\label{eq1}
M(\Gamma(|C_n|, K, R))\geqslant\delta
\end{equation}
for sufficiently large $n.$

\medskip
Let us to show that
\begin{equation}\label{eq11}
\Gamma(|C_n|, K, R)>\Gamma(g(\sigma_n), K, R)
\end{equation}
for sufficiently large $n\in {\Bbb N},$ where $\sigma_n$ denotes the
cut of $D,$ which corresponds to a domain $d_n.$ Due to the
relation~(\ref{eq3}) and item~1) of Lemma~\ref{lem1},
$C(\widetilde{f}, P)=C(f,
P)=\bigcap\limits_{k=1}^{\infty}\overline{f(d_k)}\subset
\partial D_* \,,$ therefore, $f(d_n)\cap
K^*=\varnothing$ for sufficiently large $n\in {\Bbb N}$ and any
compactum $K^*\subset D_*.$ In this case, under some $n_0\in {\Bbb
N}$ and all $n>n_0,$ we obtain that $f(d_n)\cap
H^{\,-1}(K)=\varnothing.$ Since $H$ is a homeomorphism and
$\widetilde{f}(d_n)=f(d_n)$ for such $n\in {\Bbb N},$ it follows
that $g(d_n)\cap K=\varnothing,$ $g=H\circ \widetilde{f}.$ Let now
$\gamma\in \Gamma(|C_n|, K, R).$ Since $|C_n|\subset g(d_n),$ by the
proven above $|\gamma|\cap g(d_n)\ne\varnothing\ne |\gamma|\cap
({\Bbb C}\setminus g(d_n)).$ In this case, by~\cite[Theorem~1.I.5,
$\S\,46$]{Ku}
\begin{equation}\label{eq11A}
|\gamma|\cap (\partial g(d_n)\cap R)\ne\varnothing\,.
\end{equation}
Let us now establish that
\begin{equation}\label{eq11B}
\partial g(d_n)\cap R\subset g(\sigma_n)\,.
\end{equation}
First of all, by~(\ref{eq11A}) it follows that $\partial g(d_n)\cap
R\ne\varnothing.$ Let $\zeta_0\in \partial g(d_n)\cap R.$ Then we
may find a sequence $\zeta_k\in g(d_n)$ such that
$\zeta_k\rightarrow \zeta_0$ as $k\rightarrow\infty.$ Since
$\zeta_k\in g(d_n),$ we may find $\xi_k\in d_n$ such that
$g(\xi_k)=\zeta_k.$ Since by the assumption $\overline{D}$ is a
compactum in ${\Bbb S},$ we may consider that $\xi_k$ is a
convergent sequence $\xi_k\rightarrow\xi_0\in \overline{d_n},$
$k\rightarrow\infty.$ If $\xi_0\in d_n,$ then $\zeta_0$ is an inner
point of $g(d_n)$ by the openness of $g,$ that contradicts the
choice of $\zeta_0.$ Thus $\xi_0\in \partial d_n.$ Observe also
that, $\xi_0\in D.$ Indeed, if $\xi_0\in
\partial D,$ then $f(\xi_k)=\widetilde{f}(\xi_k)$ may converge to the boundary
point of $D_*$ by the closeness of $f,$ however,
$f(\xi_k)=H^{\,-1}(g(\xi_k))=H^{\,-1}(\zeta_k)$ converges to an
inner point $H^{\,-1}(\zeta_0)\in D_*,$ because
$\zeta_k\rightarrow\zeta_0\in R$ as $k\rightarrow\infty$ and $H$ is
a homeomorphism of $U^{\,*}$ onto $R.$ The contradiction obtained
above shows that $\xi_0\in
\partial d_n\cap D,$ i.e., $\xi_0\in \sigma_n.$ Then
$g(\xi_0)=\zeta_0\in g(\sigma_n).$ Then the inclusion~(\ref{eq11B})
is established.

\medskip
Then, by~(\ref{eq11A}) it follows that $|\gamma|\cap
g(\sigma_n)\ne\varnothing.$ Therefore, the relation~(\ref{eq11}) is
also proved. By~(\ref{eq11}) and due to~\cite[Theorem~1(c)]{Fu} it
follows that $M(\Gamma(|C_n|, K, R))\leqslant M(\Gamma(g(\sigma_n),
K, R)).$ But now, by~(\ref{eq1}) is also follows that
\begin{equation}\label{eq2}
M(\Gamma(g(\sigma_n), K, R))\geqslant \delta\,,\quad n\geqslant
n_0\,,
\end{equation}
see Figure~\ref{fig1} on this occasion.
\begin{figure}[h]
\centerline{\includegraphics[scale=0.5]{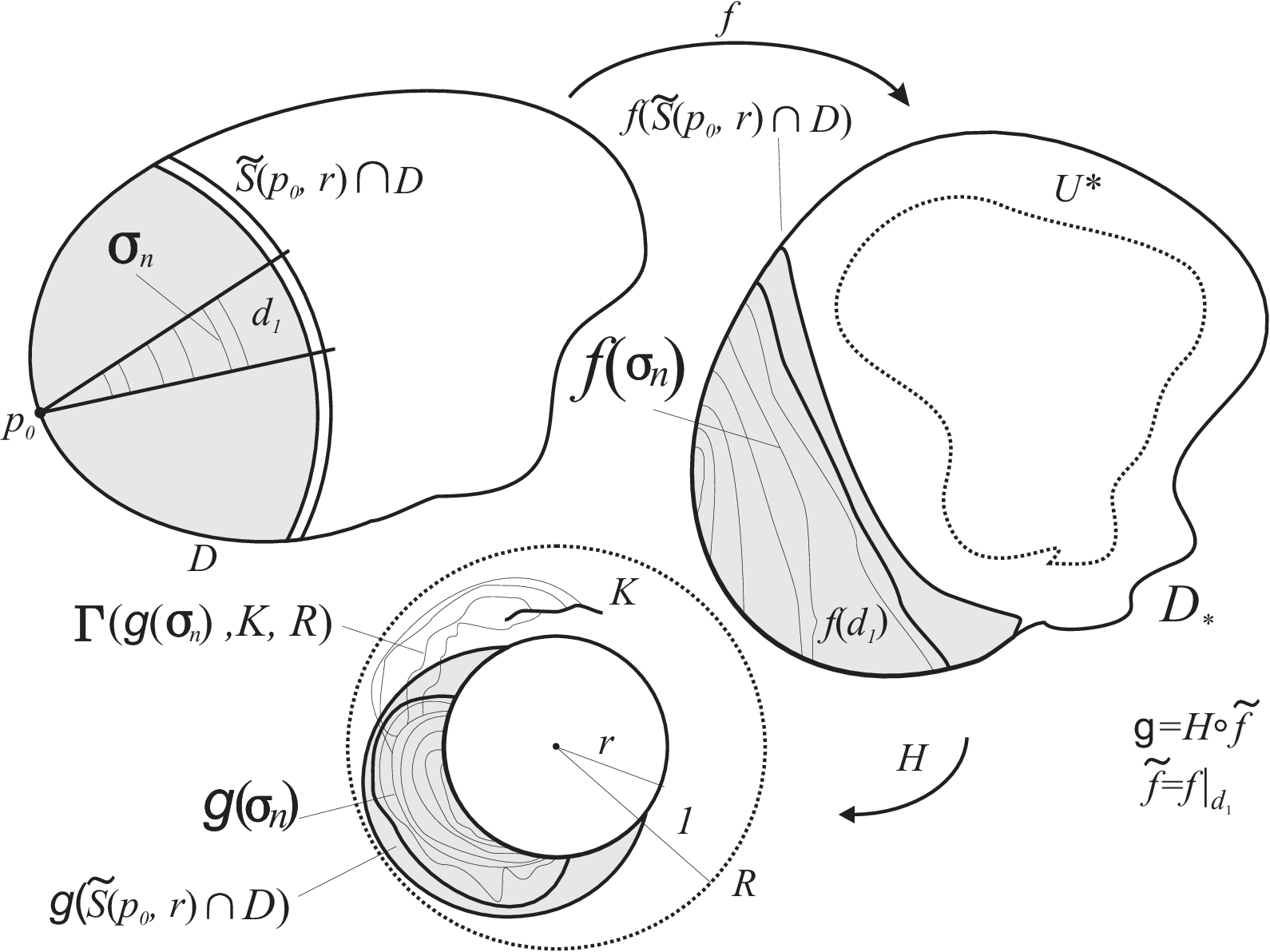}} \caption{To
the proof of Theorem~\ref{th2}}\label{fig1}
\end{figure}

\textbf{4)} Let us to show that the condition~(\ref{eq2})
contradicts with~(\ref{eq40A}). For this goal, let us to estimate
$M(\Gamma^{\,\delta}_i)$ in~(\ref{eq40A}) from above, using the
Ziemer equality and the connection between joining and separating
paths~\cite{Zi}. Let us show that, there exists $\varepsilon_1>0$
such that
\begin{equation}\label{eq9}
K\cap\overline{g(\widetilde{B}(p_0, r)\cap
d_1)}=\varnothing\quad\forall\,\,r\in (0, \varepsilon_1)\,.
\end{equation}
Assume the contrary. Then, for any $k\in {\Bbb N}$ there is $q_k\in
K\cap \overline{g(\widetilde{B}(p_0, 1/k)\cap d_1)}.$ Since $q_k\in
\overline{g(\widetilde{B}(p_0, 1/k)\cap d_1)},$ there is a sequence
$q_{kl}\in g(\widetilde{B}(p_0, 1/k)\cap d_1)$ such that
$q_{kl}\rightarrow q_k$ as $l\rightarrow\infty.$ Since $q_{kl}\in
g(\widetilde{B}(p_0, 1/k)\cap d_1),$ there is a sequence
$\zeta_{kl}\in \widetilde{B}(p_0, 1/k)\cap d_1$ such that
$g(\zeta_{kl})=q_{kl}.$ Put $k\in {\Bbb N}$ and choose $l_k>0$ such
that $|q_{kl_k}-q_k|<1/2^k.$ Without loss of generality, we may
assume that $q_k$ converges to $z_0$ as $k\rightarrow\infty.$ Then
by the triangle inequality
$$|g(\zeta_{kl_k})-z_0|\leqslant |g(\zeta_{kl_k})-q_k|+|q_k-z_0|=$$
\begin{equation}\label{eq10A}
=|q_{kl_k}-q_k|+|q_k-z_0|<1/2^k+|q_k-z_0|\rightarrow 0\,, \quad
k\rightarrow\infty\,.
\end{equation}
It follows from~(\ref{eq10A}) that $z_0\in \partial R.$ Indeed, if
$z_0\in R,$ by the continuously of $H^{\,-1}$ we obtain that
$$\widetilde{h}_*(H^{\,-1}(g(\zeta_{kl_k})),
H^{\,-1}(z_0))=\widetilde{h}_*(f(\zeta_{kl_k}),
H^{\,-1}(z_0))\rightarrow 0\,,\quad k\rightarrow\infty\,,$$
where $\widetilde{h}_*$ is a corresponding metric on ${\Bbb S}_*.$
The last contradicts with the closeness of $f$ in $D,$ because
$H^{\,-1}(z_0)\in C(f, p_0)\subset
\partial D_*$ (see~\cite[Proposition~4.3]{Sev$_1$}). At the same time,
$H^{\,-1}(z_0)$ is an inner point of $D_*$ by the assumption $z_0\in
R.$ Thus, $z_0\in \partial R,$ that contradicts with the condition
$q_k\rightarrow z_0\in
\partial R$ as $k\rightarrow\infty$ and because $q_k\in K,$ where $K$ is a compactum in $R.$ The obtained contradiction
proves~(\ref{eq9}). Then
\begin{equation}\label{eq8}
K\subset R\setminus \overline{g(\widetilde{B}(p_0, r)\cap
d_1)}\,,\quad r\in (0, \varepsilon_1)\,.
\end{equation}
In particular, (\ref{eq8}) implies that $K$ and $g(\sigma_n)$ are
disjoint for $n>n_1>n_0,$ where $n_1\in {\Bbb N}$ is such that
$r_{n_1}<\varepsilon_1.$

\medskip
\textbf{5)} Let $n>n_1.$ Observe that, for any $r\in (r_n,
\varepsilon_1),$ the set $A_r:=\partial (g(\widetilde{B}(p_0, r)\cap
d_{1}))\cap R$ separates $K$ from $g(\sigma_n)$ in $R.$ Indeed,
$$R=B_r\cup A_r\cup C_r\qquad\forall\quad
r\in (r_n, \varepsilon_1)\,,$$
where $B_r:=g(\widetilde{B}(p_0, r)\cap d_{1})$ and $C_r:=R\setminus
\overline{g(\widetilde{B}(p_0, r)\cap d_{1})}$ are open sets in $R,$
$g(\sigma_n)\subset B_r,$ $K\subset C_r$ and $A_r$ is closed in $R.$

\medskip
Let $\Sigma_n$ be a family of all sets separating $g(\sigma_n)$ from
$K$ in $R.$ Let us establish that
\begin{equation}\label{eq7A}
(\partial g(\widetilde{B}(p_0, r)\cap d_{1}))\cap R\subset
g(\widetilde{S}(p_0, r)\cap d_{1}), \quad 0<r<r_{1}\,.
\end{equation}
Indeed, let $\zeta_0\in (\partial g(\widetilde{B}(p_0, r)\cap
d_{1}))\cap R.$ Then there is a sequence $\zeta_k\in
g(\widetilde{B}(p_0, r)\cap d_{1})$ such that $\zeta_k\rightarrow
\zeta_0$ as $k\rightarrow \infty,$ where $\zeta_k=g(\xi_k),$
$\xi_k\in \widetilde{B}(p_0, r)\cap d_{1}.$ Without loss of
generality, we may assume that $\xi_k\rightarrow \xi_0$ as
$k\rightarrow\infty.$ If $\xi_0\in
\partial D,$ then by the closeness of $f$ in $D$
the sequence $f(\xi_k)$ may converge only to some boundary point
$z_1\in D_*,$ however, $H^{\,-1}(g(\xi_k))=f(\xi_k)$ converges to
some inner point of $D_*$ because $H$ is a homeomorphism,
$\zeta_k=g(\xi_k)$ and $\zeta_k\rightarrow \zeta_0\in R$ as
$k\rightarrow\infty.$ Therefore, $\xi_0\in D.$ If $\xi_0\in
\partial d_{1},$ then $\xi_0\in \sigma_{1}\subset \widetilde{S}(p_0,
r_{1}),$ that is impossible because $\xi_k\in \widetilde{B}(p_0, r)$
by the assumption, $\xi_k\rightarrow\xi_0$ as $k\rightarrow\infty$
and $r<r_{1}.$ Thus, $\xi_0\in d_{1}.$

\medskip
Now two situations are possible: 1) $\xi_0\in \widetilde{B}(p_0 ,
r)\cap d_{1}$ and 2) $\xi_0\in \widetilde{S}(p_0 , r)\cap d_{1}.$
Observe that, the case 1) is impossible because, in this case,
$g(\xi_0)=\zeta_0$ and $\zeta_0$ is an inner point of the set
$g(\widetilde{B}(p_0, r)\cap d_{1}),$ that contradicts with the
choice of $\zeta_0.$ Thus, the inclusion~(\ref{eq7A}) is
established.

\medskip
Here and below the unions of the form $\bigcup\limits_{r\in (r_1,
r_2)}
\partial g(\widetilde{B}(p_0, r)\cap d_1)\cap R$ are understood as families of
sets. Denote by $\Sigma_n$ the family of all sets which separate $K$
from $g(\sigma_n)$ in $R$ (see~\cite[section~2.3]{Zi}). Then,
by~(\ref{eq7A}), we obtain that

$$M(\Sigma_n)\geqslant
M\left(\bigcup\limits_{r\in (r_n, \varepsilon_1)}
\partial g(\widetilde{B}(p_0, r)\cap d_1)\cap R\right)\geqslant$$
\begin{equation}\label{eq5A}
\geqslant M \left(\bigcup\limits_{r\in (r_n, \varepsilon_1)}
g(\widetilde{S}(p_0, r)\cap d_1)\right)\geqslant
M\left(\bigcup\limits_{r\in (r_n, \varepsilon_1)}
g(\widetilde{S}(p_0, r)\cap d_1)\right)
\end{equation}
for $n>n_1,$ where $n_1$ is defined in item~\textbf{4)}.

By~(\ref{eq5A}) and~(\ref{eq40A}), putting $\delta=\varepsilon_1,$
we obtain that
\begin{equation}\label{eq11C}
M(\Sigma_n)\rightarrow\infty\,,\qquad n\rightarrow\infty\,.
\end{equation}
On the other hand, by the Ziemer and Hesse equalities
(see~\cite[Theorem~3.10]{Zi} and~\cite[Theorem~5.5]{Hes}), we obtain
that
\begin{equation}\label{eq6A}
M(\Sigma_n)=\frac{1}{M(\Gamma(g(\sigma_n), K, R))}\,.
\end{equation}
Now, the relations~(\ref{eq6A}) and~(\ref{eq11C})  imply that
$$M(\Gamma(g(\sigma_n), K, R))\rightarrow 0\,,\qquad n\rightarrow\infty\,,$$
that contradicts with~(\ref{eq2}). The contradiction obtained above
proves the statement of the theorem in the case, when the
relations~(\ref{eq10}) hold.

\medskip
Let us consider the case~2), namely, assume that $Q\in FMO(\partial
D).$ Let $\varphi:{\Bbb S}\rightarrow {\Bbb R},$ $\varphi(x)=0$ for
$x\not\in D,$ be a nonnegative function which has a finite mean
oscillation at a point $p_0\in \overline{D}\subset {\Bbb S}.$ By
\cite[Theorem~7.2.2]{Berd} a surface ${\Bbb S}$ is locally 2-regular
by Alhfors, so that by~\cite[Lemma~3]{Sm}
\begin{equation}\label{eq44}
\int\limits_{\varepsilon<\widetilde{h}(p, p_0)<
\widetilde{\varepsilon_0}}\frac{\varphi(p)\, d\widetilde{v}(p)}
{\left(\widetilde{h}(p, p_0)\log\frac{1}{\widetilde{h}(p,
p_0)}\right)^2} = O \left(\log\log \frac{1}{\varepsilon}\right)
\end{equation}
as $\varepsilon\rightarrow 0$ for some
$0<\widetilde{\varepsilon_0}<{\rm dist}(p_0, \partial U),$ where $U$
is some normal neighborhood of $p_0.$ Set $0\,<\,\psi(t)\,=\,\frac
{1}{\left(t\,\log{\frac1t}\right)}.$
Observe that $\psi(t)\geqslant \frac {1}{t\,\log{\frac1t}}$ for
sufficiently small $\varepsilon>0,$ therefore
$I(\varepsilon,
\widetilde{\varepsilon_0})\,:=\,\int\limits_{\varepsilon}^{\widetilde{\varepsilon_0}}\psi(t)\,dt\,\geqslant
\log{\frac{\log{\frac{1}
{\varepsilon}}}{\log{\frac{1}{\widetilde{\varepsilon_0}}}}}.$ Set
$\eta(t):=\psi(t)/I(\varepsilon, \widetilde{\varepsilon_0}).$ Them,
due to the relation~(\ref{eq44}), we may find a constant $C>0$ such
that
$$\int\limits_{\widetilde{A}(p_0, \varepsilon, \widetilde{\varepsilon_0})}
Q(p)\cdot \eta^2(\widetilde{h}(p, p_0))\
d\widetilde{v}(p)=\frac{1}{I^2(\varepsilon,
\widetilde{\varepsilon_0})}\int\limits_{\varepsilon<\widetilde{h}(p,
p_0)< \widetilde{\varepsilon_0}}\frac{Q(p)\, d\widetilde{v}(p)}
{\left(\widetilde{h}(p, p_0)\log\frac{1}{\widetilde{h}(p,
p_0)}\right)^2} \leqslant$$
\begin{equation}\label{eq46}
\leqslant C\cdot \left(\log{\frac{\log{\frac{1}
{\varepsilon}}}{\log{\frac{1}{\widetilde{\varepsilon_0}}}}}\right)^{\,-1}
\rightarrow 0\,,\quad \varepsilon\rightarrow 0\,.
\end{equation}
Then the relations~(\ref{eq10AA}) and~(\ref{eq46}) imply the
conditions~(\ref{eq45}), so that the desired conclusion follows
directly by Theorem~\ref{th2}.

\medskip
To complete the proof we may to establish the equality
$f(\overline{D}_P)=\overline{R}.$ Obviously,
$f(\overline{D}_P)\subset\overline{R}.$ Let us to show the inverse
inclusion. Let $\zeta_0\in \overline{R}.$ If $\zeta_0$ is an inner
point of $R,$ then obviously there is $\xi_0\in D$ such that
$f(\xi_0)=\zeta_0$ and, consequently, $\zeta_0\in f(D).$ Let now
$\zeta_0\in
\partial R.$ Then there is a sequence $\zeta_n\in R,$
$\zeta_n=f(\xi_n),$ $\xi_n\in D,$ such that $\zeta_n\rightarrow
\zeta_0$ as $n\rightarrow\infty.$ Since $\overline{D}_P$ is a
compactum,  we may consider that $\xi_n\rightarrow P_0,$ where $P_0$
is some prime end in $\overline{D}_P.$ Then also $\zeta_0\in
f(\overline{D}_P).$ The inclusion $\overline{R}\subset
f(\overline{D}_P)$ is proved and, therefore,
$f(\overline{D}_P)=\overline{R}.$ Theorem is proved.~$\Box$

\end{proof}

\medskip
{\it Proof of Theorem~\ref{th1}.} Observe that, $N(f, D)<\infty$
(see~\cite[Theorem~5.5]{Va$_1$}), because $f$ is an open, discrete
and closed mapping.  Note also that, an open discrete mapping
$f:D\rightarrow D_{\,*}$ with a finite distortion for which $N(f,
D)<\infty,$
$$N(f, D)\,=\,\sup\limits_{y\in{\Bbb S}_*}\,N(y, f, D)\,,$$
$$N(y, f, D)\,=\,{\rm card}\,\left\{p\in E: f(p)=y\right\}\,,$$
is a lower $Q$-mapping at any point $p_0\in \overline{D}$ for
$Q(p)=c\cdot N(f, D)\cdot K_f(p),$ where
$c>0$ is some constant depending only on $p_0$ and $D_{\,*},$ and
$K_f$ is defined in~(\ref{eq16}). In this case, the desired
conclusion follows from Theorem~\ref{th2}.~$\Box$

\medskip
{\small{\bf The datasets generated and/or analysed during the
current study are available from the corresponding author on
reasonable request.}}
\medskip


\medskip
\medskip
{\bf \noindent Evgeny Sevost'yanov} \\
{\bf 1.} Zhytomyr Ivan Franko State University,  \\
40 Bol'shaya Berdichevskaya Str., 10 008  Zhytomyr, UKRAINE \\
{\bf 2.} Institute of Applied Mathematics and Mechanics\\
of NAS of Ukraine, \\
1 Dobrovol'skogo Str., 84 100 Slavyansk,  UKRAINE\\
esevostyanov2009@gmail.com

\medskip
\noindent{{\bf Oleksandr Dovhopiatyi} \\
Zhytomyr Ivan Franko State University \\
Bol'shaya Berdichevskaya Str., 40 \\
Zhytomyr, Ukraine, 10 008 \\ e-mail: alexdov1111111@gmail.com}

\medskip
\noindent{{\bf Nataliya Ilkevych} \\
Zhytomyr Ivan Franko State University \\
Bol'shaya Berdichevskaya Str., 40 \\
Zhytomyr, Ukraine, 10 008 \\ e-mail: ilkevych1980@gmail.com }

\medskip
\noindent{{\bf Vitalina Kalenska} \\
Zhytomyr Ivan Franko State University \\
Bol'shaya Berdichevskaya Str., 40 \\
Zhytomyr, Ukraine, 10 008 \\ e-mail: vitalinakalenska@gmail.com }

\end{document}